\documentclass[12pt]{amsart}
\usepackage{latexsym,amssymb,amsmath,amsthm}
\usepackage{listings}
\usepackage{float}
\usepackage{graphicx}
\newtheorem{thm}{Theorem}
\newtheorem{lemma}[thm]{Lemma}

\newtheorem{coro}[thm]{Corollary}
\theoremstyle{definition}

 \newcommand{\R}{\mathbb{R}}  
 \def\G{\mathcal{G}} 
 \def\Det{\rm{Det}}   
  \def\sign{{\rm sign}}  
  
\def\B#1#2{{#1\choose #2}}

\title{The Euler characteristic of an even-dimensional graph}
\author{Oliver Knill}
\date{July 8, 2013}
\address{
        Department of Mathematics \\
        Harvard University \\
        Cambridge, MA, 02138
        }

\begin{document}
\begin{abstract}
We write the Euler characteristic $\chi(G)$ of a four dimensional
finite simple geometric graph $G=(V,E)$ in terms of 
the Euler characteristic $\chi(G(\omega))$ of two-dimensional geometric subgraphs $G(\omega)$. 
The Euler curvature $K(x)$ of a four dimensional graph satisfying the Gauss-Bonnet relation
$\sum_{x \in V} K(x) = \chi(G)$ can so be rewritten as an average $1-{\rm E}[K(x,f)]/2$ over a collection
two dimensional ``sectional graph curvatures" $K(x,f)$ through $x$. Since scalar curvature, the average of all 
these two dimensional curvatures through a point, is the integrand of the Hilbert action,
the integer $2-2 \chi(G)$ becomes an integral-geometrically defined relative of the Hilbert action functional.
The result has an interpretation in the continuum for compact $4$-manifolds $M$:
the Euler curvature $K(x)$, the integrand in the classical Gauss-Bonnet-Chern
theorem, can be seen as an average over a probability space $\Omega$ of $1-K(x,\omega)/2$ 
with curvatures $K(x,\omega)$ of compact $2$-manifolds $M(\omega)$. Also here, the Euler characteristic
has an interpretation of an exotic Hilbert action, in which sectional curvatures are replaced
by surface curvatures of integral geometrically defined random two-dimensional sub-manifolds 
$M(\omega)$ of $M$. 
\end{abstract} 
\maketitle

This is an informal note explaining a comment which slipped into \cite{IsospectralDirac2}.
It uses the observation of \cite{indexformula} that the symmetric index $j_f(x) = (i_f(x) + i_{-f}(x))/2$
at a critical point $x$ of a function has a topological interpretation as the genus $1-\chi(B_f(x))/2$ 
of a lower dimensional space $B_f(x)$. The index $i_f(x) = 1-\chi(S^-_f(x))$ is a discretisation of the
index of a gradient vector field $\nabla f$ at a critical point which by Poincar\'e-Hopf add up to 
the Euler characteristic $\chi(G)$ of $G$.
For four dimensional spaces, $j_f(x)$ is the genus of a two-dimensional compact surface $B_f(x)$ 
obtained by intersecting a small sphere $S(x)$ with the level surface of $f$ at $x$. 
Since genus is additive, we can glue the local critical surfaces $B_f(x)$ together and get for every function $f$ a 
two-dimensional graph $G(f)$ whose genus is the sum of the indices.  Poincar\'e-Hopf assures that the Euler 
characteristic of this surface is related to the Euler characteristic of the entire space. If we integrate over a probability space of
functions $f$, the expectation of the curvatures $K(x,f)$ of these surfaces at a vertex $x$ is related to the Euler curvature
$K(x)$. Because scalar curvature classically is an average over all sectional curvatures, this brings Euler characteristic
in the vicinity of the Hilbert-Einstein action in differential geometry and suggests to 
search for graphs which maximize or minimize the Euler characteristic. The question is then whether there are local rules
similar than the vacuum Einstein equations which assure that the Euler characteristic is extremal
and what are the geometric properties of the extrema. While we can not yet answer this  yet, we will comment on 
it anyway. In the graph case, where traditional tensor calculus is absent, it is natural to look at 
the Einstein tensor $T(v,e) = R(e) - S(v)$ which involves the Ricci curvature $R$, the average of wheel 
graph curvatures through an edge $e$ and scalar curvature $S(v)$ the average of wheel graph curvatures through 
a vertex $v$. An icosahedron for example satisfies $T(v,e)=0$ for all vertices $v$ and edges $e$
the reason being that $R(e)=1/6$ and $S(e)=1/6$ everywhere. Having $T$ zero everywhere, it is an Einstein graph. 
The just mentioned notions for Ricci and scalar curvature for graphs are rather rigid and can not be deformed by quantum 
dynamics or unitary symmetries like isospectral Dirac deformations. We are going to replace them therefore. 
The starting point is that the Euler curvature $K(x)$ can be deformed because it is the expectation of indices $i_f(x)$. 
We can define now two-dimensional curvatures $K(x,f)$, replacing sectional curvatures and get in a familiar way 
Ricci curvature and so an Einstein tensor as an expectation over $f$. The Einstein equations then define the mass 
tensor for any even-dimensional geometric graph.  The last point is important;
by defining the mass distribution from the geometry data, we assure that geometry remains the only input. 
The hope is that the Euler characteristic as a variational
problem selects interesting geometries with interesting mass distributions. There is a second aspect we can
address: the geometric ideas also allow for four dimensional graphs to define a notion
of an action given as the expectation value of the genus of a surface $B_f(\gamma)$ connecting two vertices. Also this
notion is deformable and can be used to deform the geodesic distance which in general a much too small radius of 
injectivity for discrete networks. But lets start from the beginning: \\ 

Inductively, a finite simple graph $G=(V,E)$ is called geometric of dimension $d$, if every unit 
sphere $S(x)$ is a $(d-1)$-dimensional geometric graph of Euler characteristic 
$\chi(S(x)) = 1-(-1)^d$. The induction assumption is that any zero-dimensional graph - a 
graph without edges - is geometric. One could strengthen the assumption and ask that each $S(x)$ is a Reeb sphere,
a $(d-1)$-dimensional graph which admits an injective function with exactly two critical points of index $1$;  but we do not do that
because it is not needed. Examples of geometric graphs are sufficiently nice triangularizations of $d$-dimensional smooth manifolds. 
Given a real-valued injective function 
$f:V \to \R$ and $c \in \R$ different from any of the values of $f$, we can partition the vertex set $V$ into two sets 
$V^+_f =  \{ x \; | \; f(x)>c \; \}$ and $V^-_f= \{ x \; | \; f(x)<c \; \}$.
Define the hyper-surface graph $G_f$ whose vertices are edges in $E$ on which $f(y)-f(x)$
changes sign and whose edges are triangles in $G$ on which $f(y)-f(x)$ takes different signs.
We think about $G_f$ as the discrete analogue of the level surface $\{ f = c \; \}$ contained inside the graph and passing
through $x$.  We have proved in \cite{indexformula} that the graph $S(x)_f$ 
is a polytop: it can be completed in a canonical way to become a $(d-1)$-dimensional geometric graph $B_f(x)$ 
and have seen that, in general, the index formula  $j_f(x)= 1-\chi(S(x))/2- \chi(B_f(x))/2$ holds,
where $j_f(x) = (i_f(x)+i_{-f}(x))/2$ and $i_f(x) = 1-\chi(S_f^-(x))$. 
For odd-dimensional geometric graphs, this gives $j_f(x) = -\chi(B_f(x))/2$ which leads to the statement 
that curvature at $x$ - the expectation of $j_f(x)$ over all functions $f$ - 
is zero everywhere for odd-dimensional graphs. For even dimensional geometric graphs, 
the formula becomes $j_f(x) = 1-\chi(B_f(x))/2$ which is a ``genus" with an additivity property.
For four-dimensional graphs in particular, we can write $\chi(G)$ in terms of an expectation of Euler characteristics of 
two dimensional graphs. For six dimensional graphs, since it reduces to the expectation of Euler characteristic 
of four dimensional graphs, which each can be reduced to a sum of two dimensional graphs, we can again write 
$\chi(G)$ as an expectation of two-dimensional graphs, etc. Not to complicate things, we stick mainly 
to four dimensions but inductively, we can reduce every even dimensional graph to two dimensional graphs.
The genus can be spotted in \cite{Spivak5} on page 234 in a two-dimensional case: 
the relation between index $i$ and order $s$ of a saddle is $s=1-2i$ or $i=1-s/2$, where generically the intersection
of the level curve through the critical point intersected with a small sphere is $4$. \\

The Euler characteristic $\chi(G) = \sum_{k=0} (-1)^k v_k$ of a graph $G$ is a super counting function which satisfies 
$\chi(G \cup H) = \chi(G) + \chi(H) - \chi(G \cap H)$, like cardinality. Indeed, the
identity is the exclusion-inclusion picture applied in parallel to the number $v_k$ of all sub-simplices of 
dimension $k$. It implies for example that the number $\chi(G)-1$ is additive when merging two graphs along a simple 
vertex. It also assures that if two geometric 2-dimensional graphs are glued along a contractible circle 
and the discs bounded by the circles are taken out on both sides, then the genus $g=1-\chi(G)/2$ is 
additive as long as we apply it to two-dimensional geometric graphs which are surfaces.  Similarly, if 
two geometric $4$-dimensional graphs are glued along a 
$3$-dimensional manifold which bounds contractible pieces on both sides, then $1-\chi(G)/2$ is additive. For example,
if two $4$-spheres of Euler characteristic $2$ are glued along a $3$ sphere which bounds contractible sets in both spheres,
we obtain a larger sphere of Euler characteristic $2$. If a $2$-torus is glued along a contractible circle to a given surface, 
then the genus, the number of holes increases by $1$ because one more hole is added. 
We see that the index $j_f(x)$ of a vertex of a four-dimensional geometric graph
is equal to the genus of the hyper-surface polytop $B_f(x)$ in the unit sphere, which is a two-dimensional 
graph. The additivity allows us to glue the individual graphs together, as long as we glue only disjoint graphs.
This produces for every tree $t$ connecting all the vertices in the graph a single two-dimensional subgraph $G(f,t)$.
By the Kirkhoff matrix tree theorem, there are $\Det(L_0)/n$ trees, where $L_0$ is the 
scalar Laplacian and $\Det(A)$ is the product of nonzero eigenvalues of $A$.
We could also consider $\chi(G) \Det(L_0)/n$ as a functional for graphs
because it is a sum over all possible ``paths", we can also just average over all possible maximal trees and 
see $\chi(G)$ itself as fundamental also because $\Det(L_0)$ is not invariant under homotopy deformations.

\begin{lemma}[Genus lemma]
Both for $4$-manifolds and geometric $4$-graphs, the symmetric Morse index $j_f(x)$ at a critical point $x$
of $f$ is equal to the genus of the $2$-manifold or geometric $2$-graph $B_f(x)$ defined by intersecting
the level surface through the critical point with a sphere. 
\end{lemma}

By index expectation and Gauss-Bonnet, we have a geometric interpretation of Euler curvature, the integrand of
Gauss-Bonnet-Chern: 

\begin{coro}
The Euler curvature $K(x)$ at a point $x$ is the genus expectation for random surfaces obtained at $x$. 
\end{coro}

As explained, we understand $S(x)=S_r(x)$ as the geodesic sphere of sufficiently small radius if
we are in the Riemannian manifold case and $\{y \; | f(y) = f(x) \; \}$ to be a subgraph of
the simplex graph $\G$ of $G$ if we are in the graph case. Still in the graph case, 
we understand that $B_f(x)$ always has been completed and hence has become geometric. It is not yet clear how much the 
above lemma can be generalized to general finite simple graphs. The genus becomes in general half an integer as the
triangular graph shows already. The graphs $B_f(x)$ can still be defined, but it is not yet clear how to complete them
nor to glue the graphs $B_f(x)$ from various critical points in an additive way. 
It would be interesting to have that because it would express
the Euler characteristic of a finite simple graph in terms of the average of the unit sphere Euler characteristic
and a graph $B_f$ of smaller dimension. In practical terms, it is not really necessary to have a geometric interpretation: 
Poincar\'e-Hopf already reduces the computation of the Euler characteristic of $G$ to the computation of Euler characteristics of 
subgraphs of the unit spheres and this is by far the fastest way to compute $\chi(G)$. 
Done on a computer,  it beats every other method at great length. It essentially makes
the computation of Euler characteristic a polynomial task from a practical point of view.
(There are graphs with high dimension,
where unit spheres are large and where the inductive computation using spheres does not break the task 
down quick enough so that it is not polynomial in general), while other methods are exponential. \\

In the graph case, there is a natural probability space of scalar functions on the vertices: take functions 
which take random uniformly distributed values in $[-1,1]$ on each vertex and for which the values at different vertices are
independent. Since the sum $\sum_{x \in V} j_f(x)$ is always the Euler characteristic by Poincar\'e-Hopf, 
we can interpret the Poincar\'e-Hopf theorem for $4$-dimensional graphs $G$ as the fact that the Euler characteristic of $G$
is equal to the Euler characteristic of any choice of a two dimensional graph $G(f)$. The later can be thought of as 
a ``string", a two-dimensional surface in the four dimensional ``space-time". 
This remains true for the expectation value of $\chi(G(f))$ over a probability space of functions. 
We have shown that this is equal to the curvature. And it remains true if we average
over a probability space of trees $t$. Now, look at a vertex $x$ and consider the curvatures $K(x,\omega)$ 
of all the two-dimensional graphs $G(\omega)$ passing through an edge through $x$. 
We consider any of the $K(x,\omega)$ as a choice of a ``sectional curvature". 
We have shown that that the average $K(x) = 1-{\rm E}[K(x,\omega)]/2$ is the Euler curvature. In other words,
we have conceptionally placed the Euler curvature $K(x)$ in the vicinity of scalar curvature 
and Euler characteristic in the vicinity of the Hilbert-Einstein action. 

\begin{thm}[Euler characteristic is Hilbert-Einstein]
For geometric $4$-dimensional graphs $G$, the Euler characteristic of $G$
is equal to the Euler characteristic of each of the embedded $2$-dimensional random geometric
graphs $G(\omega)$. The Euler curvature $K(x)$ at a vertex $x$ is the expectation of the 
curvature expressions $1-K(x,\omega)/2$ of the random two dimensional graph $G(\omega)$ at $x$.  
\end{thm} 

Of course, this match is only conceptional since the curvatures $K(x,\omega)$ are not sectional curvatures 
and there will hardly be closer link as the classical Hilbert action is a real number while Euler 
characteristic is an integer. Also, the Hilbert action is not a homotopy invariant, while Euler characteristic 
is.  The statement however should add weight to the believe that Euler characteristic is an important functional 
in the graph case and in the manifold case under some curvature and volume constraints.  \\

Lets look at the second variational problem in relativity, the search for geodesics, when geometry 
is fixed. Also this functional can be modeled over a  probability space chosen on functions and so become deformable:
the Euler characteristic of a two-dimensional surface defines a functional for the set of graphs 
$\gamma: x_0,x_1, \dots ,x_n$ connecting two vertices $a,b$ in a four-dimensional geometric graph $G$. If we glue the polytopes 
$B_f(x_j)$ along the path, we get a two-dimensional graph $G(f,\gamma)$. The expectation of the Euler characteristic 
$\chi(G(f,\gamma))$ when averaging over all functions $f$ gives an action $S(\gamma)$. It is not necessarily an integer after 
averaging and so more flexible than the usual geodesic distance in a graph. If $G$ is a $3$-dimensional geometric graph, 
then the graphs $B_f(x_j)$ are one dimensional; gluing them together produces then a one-dimensional graph with several components. But the 
Euler characteristic is always zero. We can now look at the path which minimizes the action $|\gamma| - \epsilon S(\gamma)$, which
is a metric for small enough $\epsilon>0$. One could also look at path integrals $\exp(i S(\gamma))$ over all possible paths $\gamma$.
While arc-length does not deform and has a rather small radius of injectivity, the new metric changes, when deforming 
the Dirac operator on the graph. Why are we unhappy about the usual geodesic metric on a graph? 
Whenever a graph has two triangles sharing a common edge, then there are already two geodesics of length $2$ connecting vertices
in this kite graph: the caustic is close. But we would like to have similar differential geometry than in the continuum. 
For a triangularizations of a sphere like a two dimensional Buckminster type graph, we would like to  to have the 
caustics appear near the antipode. In short: we want a nicer and more flexible exponential map on a graph which is even more
sensitive to curvature. The genus action is a real number which can now distinguish the shortest connection. 
It plays well with curvature because negative curvature will produce surfaces $G(f,\gamma)$ with large genus. 
We have hopes that the new tool also allows to prove things better like classical results in differential geometry, 
for example Hadamards theorem for graphs with negative curvature, or other theorems where the exponential map 
plays an important role. \\

Anyway, we see that both pillars of general relativity, the task to ``get geometry from matter" or the task to see 
``how matter moves in a given geometry" can be framed within graph theory in such a way
that a unitary deformation on the space of functions deforms both the geometry as well as the exponential map respectively
the geodesic flow. The integral-geometric point of view adds more flexibility in the discrete, even without 
the availability of tensor calculus. (Integral geometry of course is rooted deeply in differential geometry, not at least by the
influence of Blaschke and through Chern). 
We want flexibility because interesting in geometry is done by deformation, Ricci deformation is only the latest 
example of how one can see that deformation is a powerful variant of induction or decent.
Having found an integrable deformation in Riemannian geometry \cite{IsospectralDirac2}, we of course hope that 
this might become useful. In any case, the notions considered here go well with such deformations.  \\

So far, we have looked mainly at four dimensional graphs. 
The action $S(\gamma)$ are even useful for two-dimensional graphs, where
$B_f(x)$ is a zero-dimensional graph and $j_f(x) = 1-\chi(B_f(x))/2$ gets larger if the curvature is getting smaller.
We can look therefore at metrics $n- c \sum_{k=0}^n j_f(x_k)$, where $c$ is sufficiently small, to have a metric. Now, the distances
are made larger at places with negative curvature. Again, the distance has become more flexible. 
For any measure on the space of functions, we get a distance. Most of these measures will now produce metrics for which 
the radius of injectivity is larger. One could even use this to select out measures: find a  measure on functions 
such that the sum of the radii of injectivity is maximal. For an icosahedron for example, we want the radius of 
injectivity to be 3 for every vertex so that wave fronts only focus at the antipode. 
Rather than artificially weight the graph by changing the lengths of the edges, the change is made
which is more in line with curvature and therefore natural. \\

For Riemannian manifolds $M$, there is a similar story. Again, we can replace tensor analysis with an integral geometric
framework. We can still show for $4$-manifolds that the Euler curvature $K(x)$ - the integrand of the Gauss-Bonnet-Chern 
theorem - is the expectation $1-K(x,\omega)/2$ involving  curvatures $K(x,\omega)$ of two-dimensional surfaces in $M$. 
This is indeed true, even so we have not yet found an intrinsic and natural probability space on the space of scalar functions.
We have experimented (*) with different probability spaces of Morse functions. 
One possibility is to take the manifold $M$ with normalized volume measure as
the probability space and take for every point $x \in M$ the heat kernel function $f(y) = [e^{-\tau L_0}]_{xy}$ then possibly
integrate over a probability measure of $\tau's$ and the volume measure for $x$. But we have not yet verified which measure
leads to curvature as an expectation. Work like \cite{Banchoff67} suggests an other approach: 
Nash embed the Riemannian manifold $M$ into an ambient 
Euclidean space $E$ and look at all linear functions on $E$ with unit gradient. The so induced functions on $M$ produce a finite
dimensional probability space. The embedding approach is elegant but is not intrinsic yet. \\

\begin{tiny}
(*) Before discovering the link between index and curvature in the discrete, we have experimented 
numerically with heat kernel approaches in the continuum which gives Morse functions for each $x$, 
where $M$ itself is the probability space. This is not easy to explore numerically in the 
continuum in Riemannian setups, because many geodesics have to be computed to estimate the heat kernel
$K(t,x,y)=\sum_n e^{-\lambda_n t} f_n(x) f_n(y) = \exp(-t L)(x,y)$, the fundamental solution of the heat 
equation $(d_t + L) f=0$. 
While it is likely that the index density of the heat kernel is Euler curvature, is still unproven. The experimental evidence 
is not conclusive and I myself got distracted by graph theory. 
The intuition is that the {\bf diffusion distance} $d_t(x,y)^2 =K(t,x,x)+K(t,y,y)-2K(t,x,y)$
is a "quantum distance" between two points. Unlike the geodesic distance, it is smooth everywhere on $M$
and gives in the limit $t \o 0$ the usual distance by Varadhan's lemma $\lim_{t \to 0} t \log(K(t,x,y)) = -1/2 d(x,y)^2$.
In the Euclidean case, $K(t,x,y)=(4\pi t)^{-n/2} \exp(-(x-y)^2/(4t))$.
Because curvature can be recovered from the heat kernel by $K(t,x,x) = (4\pi t)^{-n/2}(1+K(x)/6 t+ O(t^2) ...)$
we expect critical points of the distance function to be near points where curvature is large with positive index near
maxima or minima and saddle points near points where curvature is negative.
For fixed $x$, the {\bf heat kernel signature} function $f_y(t): t \to K(t,x,y)$ is called a {\bf heat kernel map}.
We have numerically constructed the heat kernel $K(t,x,y)$ by running Brownian
paths from each point x for some time $t$ and look at the density of the end points. Unfortunately this does not
give an accurate picture, because we have to find the critical points in each case and then run
things from many initial points. Brownian motion on a Riemannian manifold is a Markov process whose transition density function is the
heat kernel associated with the Laplace-Beltrami operator $L$. On a compact manifold the flow 
is stochastically complete and satisfies a Dynkin formula ${\rm E}[f(X_t)] = f(x) + {\rm E}[\int_0^t L f(X_s) \; ds]$.
Intuitively the heat approach makes sense: if we take a small bump with positive curvature on a manifold, then the level curves of
the Green function will have points of positive index near the bump and points of negative index near
the rim where curvature is negative. 
The index density $I_t(x)$ of the heat kernel might a priory depend on $t$. By Poincar\'e-Hopf
we are always led to a ``curvature function" which gives when integrated the Euler characteristic. The question is whether
it is the traditional curvature. 
One question we would like to answer first is: if we change a manifold outside a neighborhood of a point $x$, does the 
index density change near that point?
If not, this would add confidence to a conjecture that the average index density of 
all the heat signature functions is equal to the Euler curvature for all $t>0$. 
\end{tiny}
\vspace{4mm}

Lets take a $4$-manifold and let $S(x)$ denote sufficiently small geodesic sphere $S_r(x)$.  
We get two-dimensional manifolds $B_f(x)$ and have the classical symmetric index $j_f(x) = (i_f(x) + i_{-f}(x))/2$ for 
$4$ manifolds written as the genus $1-\chi(B_f(x))/2$. If $B_f(x)=\emptyset$ like for maxima or minima of $f$, the genus is $1$. 
Taking the probability space for granted, we can for every $f$ take a sufficiently fine triangularization $T$ of $M$ which is
a graph of the same dimension and for which all critical points belong to vertices. Then chose a 
tree $t$ in $T$ which contains these points. We can now glue a two-dimensional surface $M(f,t)$.
Denoting the elements of the probability space $\Omega$ with $\omega = (f,t(f))$ of $\Omega$, we can 
now define curvatures $K(x,\omega)$ of the $2$-manifold $M(\omega)$ through $x$. The expectation of $1-K(x,\omega)/2$ 
produces the Euler curvature $K(x)$. We see that also here, the Euler curvature is an average over 
sectional curvatures and the Euler characteristic of $M$ the expectation over the Euler characteristics over a 
class of two-dimensional sub-manifolds in $M$. 
The upshot is that for any probability measure on the space of Morse functions and any choice of 
triangularizations $t(f)$ for each function, there is a curvature $K(x)$ which integrats to $\chi(M)$ such that $K(x)$ 
is an expectation of curvatures of two dimensional surfaces passing through $x$. 
This shows that also in the continuum in even dimensions, we can see the Euler characteristic as an exotic Hilbert functional. \\

Besides the fact that a natural intrinsic probability space still needs to be found, also the gluing procedure of different $B_f(x,\omega)$ 
is not canonical. The curvatures $K(x,\omega)$ are not a sectional curvature because the surface curvature of the 2-manifolds is different 
from the sectional curvature in the tangent direction to the 2-manifold. Lets look at the homogeneous 4 sphere $M$ embedded 
in $E=R^5$. Linear functions induce Morse functions on $M$ which have a maximum and minimum. At both points $B_f(x)$ is
empty and the genus is $1$. At all other points $B_f(x)$ is a two sphere of genus $0$. Gluing them all together along a tree
of a triangularization produces one big 2 sphere $M(\omega)$ of genus $0$. What contributes to the curvature is the
empty space near the extrema. We can visualize $M(\omega)$ having flat parts there so that $1-K(x,\omega)/2=1$ there.\\

\begin{figure}
\parbox{15.4cm}{
\parbox{7cm}{\scalebox{0.22}{\includegraphics{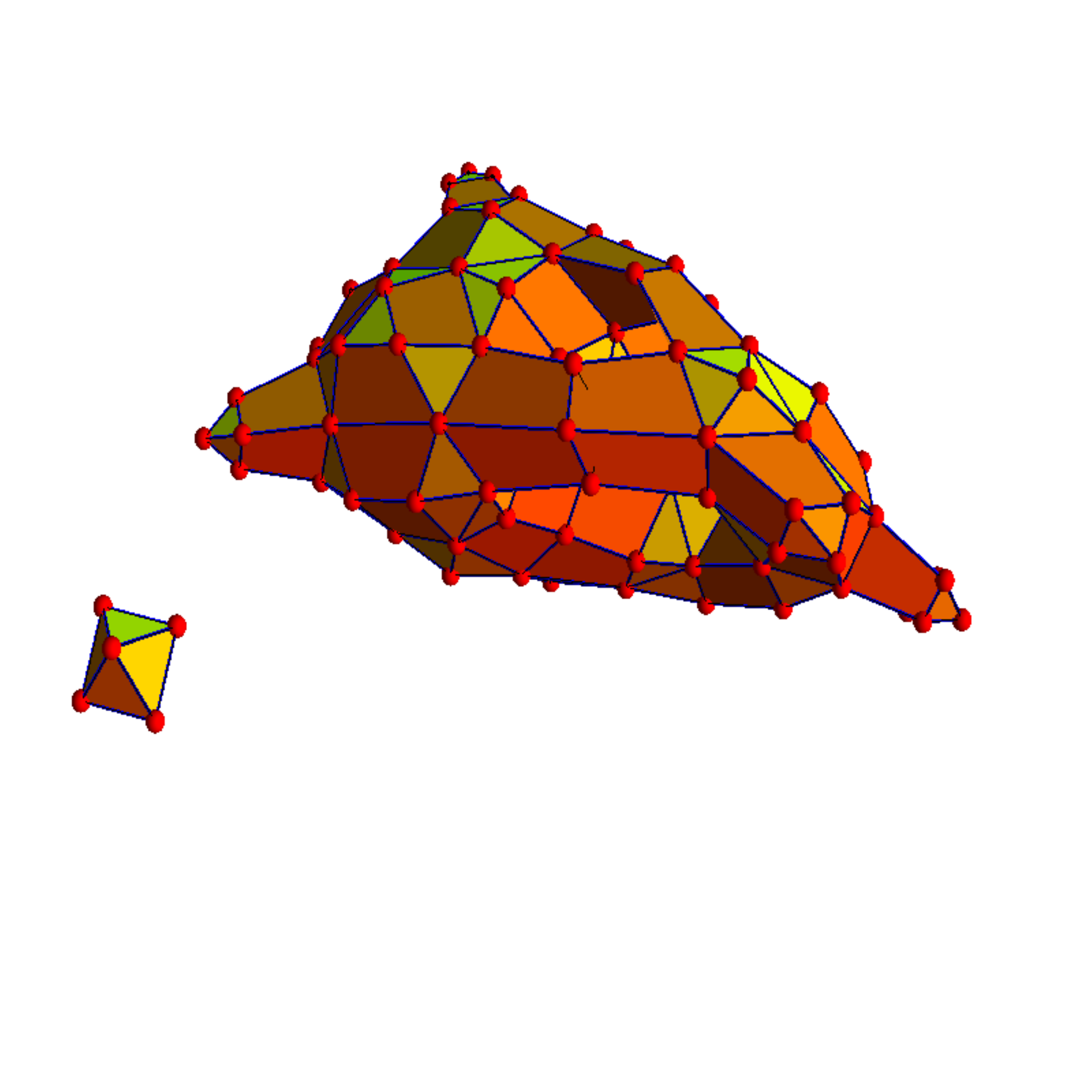}}}
\parbox{7cm}{\scalebox{0.22}{\includegraphics{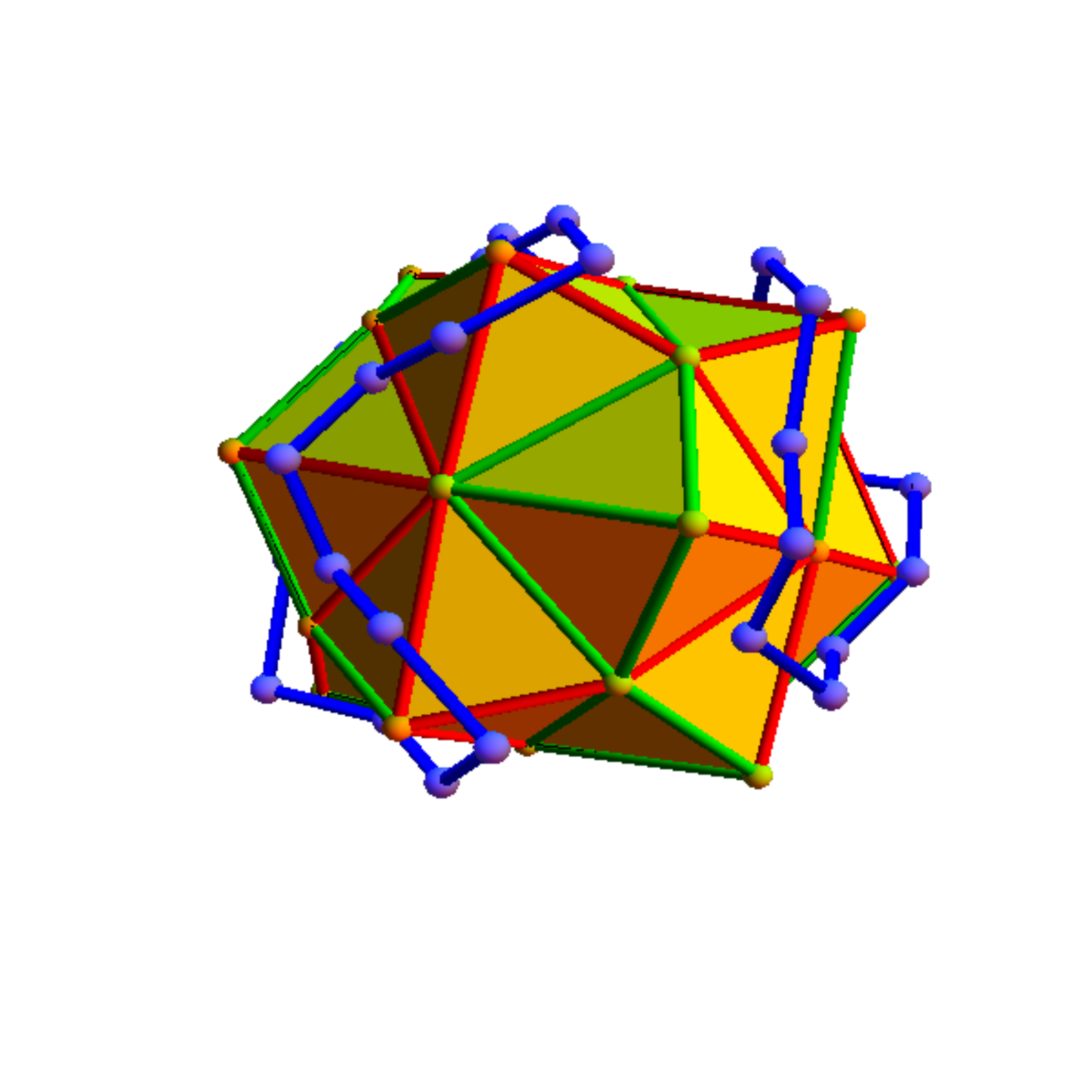}}}
}
\caption{
An example of a random graph $G_f(x)$, a polytop. We are in the situation where $S(x)$ is 
a three dimensional sphere. We chose a random function $f$ on the vertices and computed
the two dimensional surface. It does not have to be connected. In this case, one is
a sphere, the other has higher genus.  The right figure shows the case of a three dimensional $G$,
where the unit sphere $S(x)$ is two dimensional. Then $B_f(x)$ is one
dimensional graph. Also this does not need to be a single closed path.  }
\label{wk}
\end{figure}

What happens in odd dimensions? 
For odd dimensional graphs or manifolds, the reduction for any function $f$ ends up with 
one-dimensional closed graphs $B_f(x)$ which have Euler characteristic $0$. This has been used to
prove that the curvature for odd dimensional geometric graphs is always constant zero. 
Euler curvature for an odd dimensional manifold $M$ is usually not defined but could be by setting it 
to be constant $0$. Lets look at a three dimensional manifold, a function $f$ and a tree $t$ of a triangularization.
The one dimensional manifolds $B_f(x)$ can as before be glued together to form a collection of closed loops. 
Because $B_f(x)$ is a collection of loops on $S(x)$, a two dimensional sphere,
they are not knotted in the ambient space. It is again true that curvature is an
expectation of Euler curvatures, the statement is just trivial now and we do 
not get an interesting identity.  \\

So, also in the manifold case, we have gained generality by writing Euler curvature $K(x)$ as an expectation of curvature expressions 
$1-K(x,\omega)/2$ of smaller dimensional manifolds. The probabilistic setup makes sense also for manifolds $M$ 
which are no more smooth and since Euler characteristic is a robust homotopy invariant, things are pretty deformable. 
Given a homeomorphic deformation of $M$, the curvature can be pushed along simply by pushing forward the probability measure
on functions. 
Gauss-Bonnet-Chern can so be generalized to polytopes (where it is of course well known) 
or even more singular objects. We can for example deform 
a manifold to become a piecewise linear manifold for which curvature is located on the vertices. The setup is not only 
robust under deformations, it is even robust under homotopies. We can for example define a homotopy
which changes a manifold $M$ to $M \times [0,1]$. This thickened manifold has higher dimension and a boundary 
but its Euler characteristic is the same. Of course, also the two dimensional surfaces will be thickened and become
three dimensional manifolds with boundary of the same Euler characteristic. The probabilistic picture allows to 
push ideas of curvature and results like Gauss-Bonnet-Chern into areas, where tensor analysis is no more available.
We can even forget about the Euclidean fillings in the polytopes and end up with the graph case. Thats what topologists
have done since the very beginning. \\

Now the race is on to find local ``Sarumpaet rules" \cite{Egan} which play the role of the Einstein equations in the continuum 
and which are true that the Euler characteristic is extremal. Because Euler curvature is similar to 
scalar curvature, we have to replace the ingredients of the Einstein vacuum equations $R - g S = 0$.
The obvious step is to assume Ricci curvature $R(e)$ to be the average over all sectional
curvatures of two dimensional surfaces which contains $e$ and $S$ is the scalar curvature, the average of all $R(e)$
with $e$ connected to $v$. Lets call a rule local at a point $x$
if it affects only properties $p(y)$ for $y$ in the unit ball $B(x) = \{x \} \cup S(x)$ of the graph and properties $p(y)$
are local in the sense that they are same at $y$ if $G$ is replaced by $B(y)$. 
In other words, a rule is local if applied at $y$
is the same if the graph is replaced by the ball $B_2(y)$ of radius $2$. \\

{\bf The Sarumpaet Problem:} are there local rules which 
are satisfied by every graph of order $n$ which has extremal $\chi(G)$ among all other graphs 
of order $n,n+1$. Are they related to Einstein type equations? \\

The corresponding question for manifolds is tricker because Euler characteristic is an integer
which does not change under topological or even homotopy deformations. The problem also only makes 
sense for even dimensional manifolds because $\chi(M)=0$ for odd dimensional ones. 
Allowing to rip apart a manifold 
without giving a notion of what ``local perturbation" means would not make sense. Also, for manifolds, there
are no extrema without bounding something like curvature or volume:
for $4$-manifolds, one has $\chi(M) = 2+b_2-2b_1$ by Poincar\'e duality which shows that $\chi$ is unbounded also
above: a 4 dimensional Swiss cheese with lots of holes has large Euler characteristic. 
One could ask for manifolds with maximal or minimal Euler characteristic among all manifolds of dimension 
$d$, fixed volume and for which all sectional curvatures are bounded in some interval.
For manifolds, we might collapse the manifold to a nice triangularization first which has the same Euler characteristic
and then work with variations of the graph. The questions for graphs are definitely easier and more natural. 

\begin{figure}
\parbox{15.4cm}{
\parbox{7cm}{\scalebox{0.22}{\includegraphics{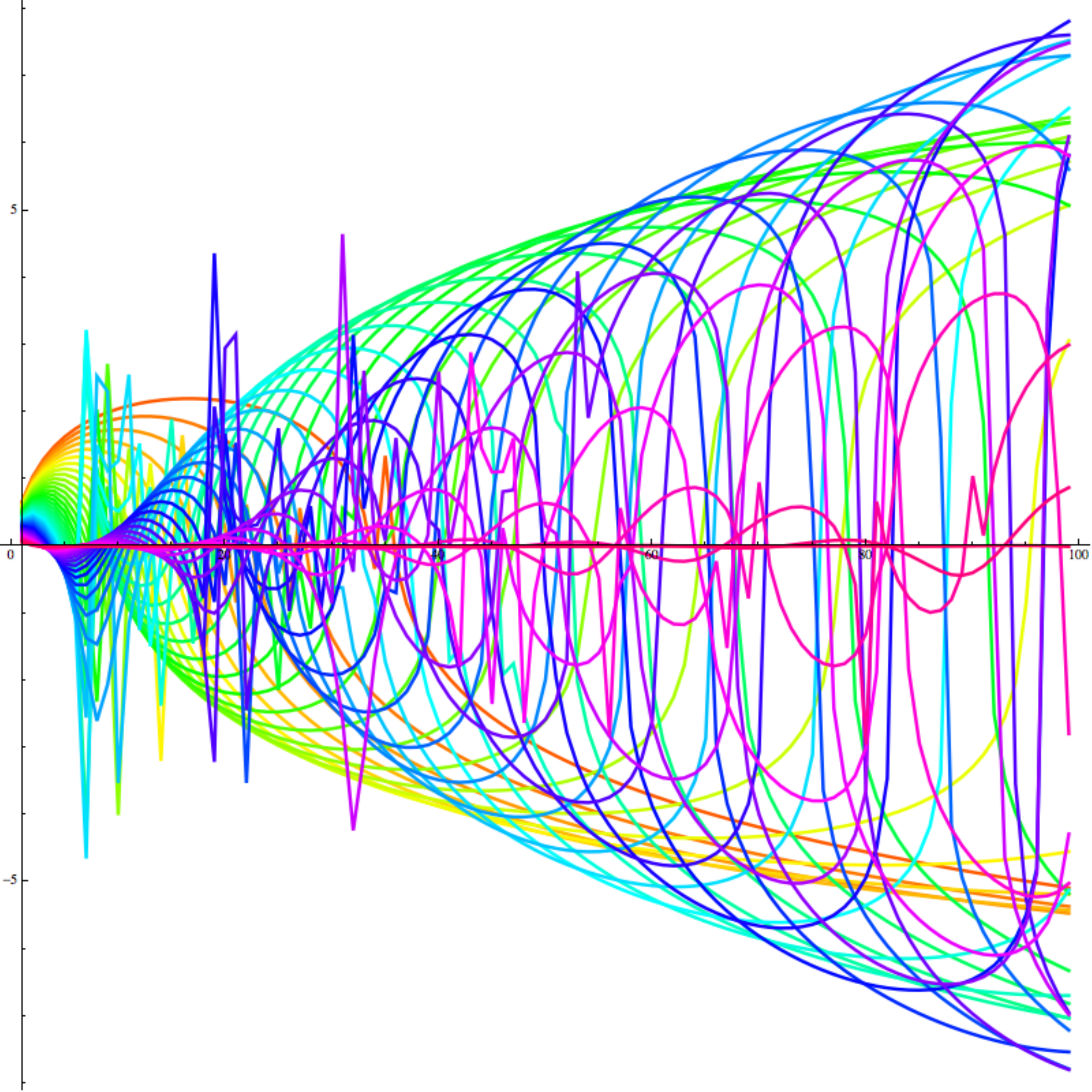}}}
\parbox{7cm}{\scalebox{0.22}{\includegraphics{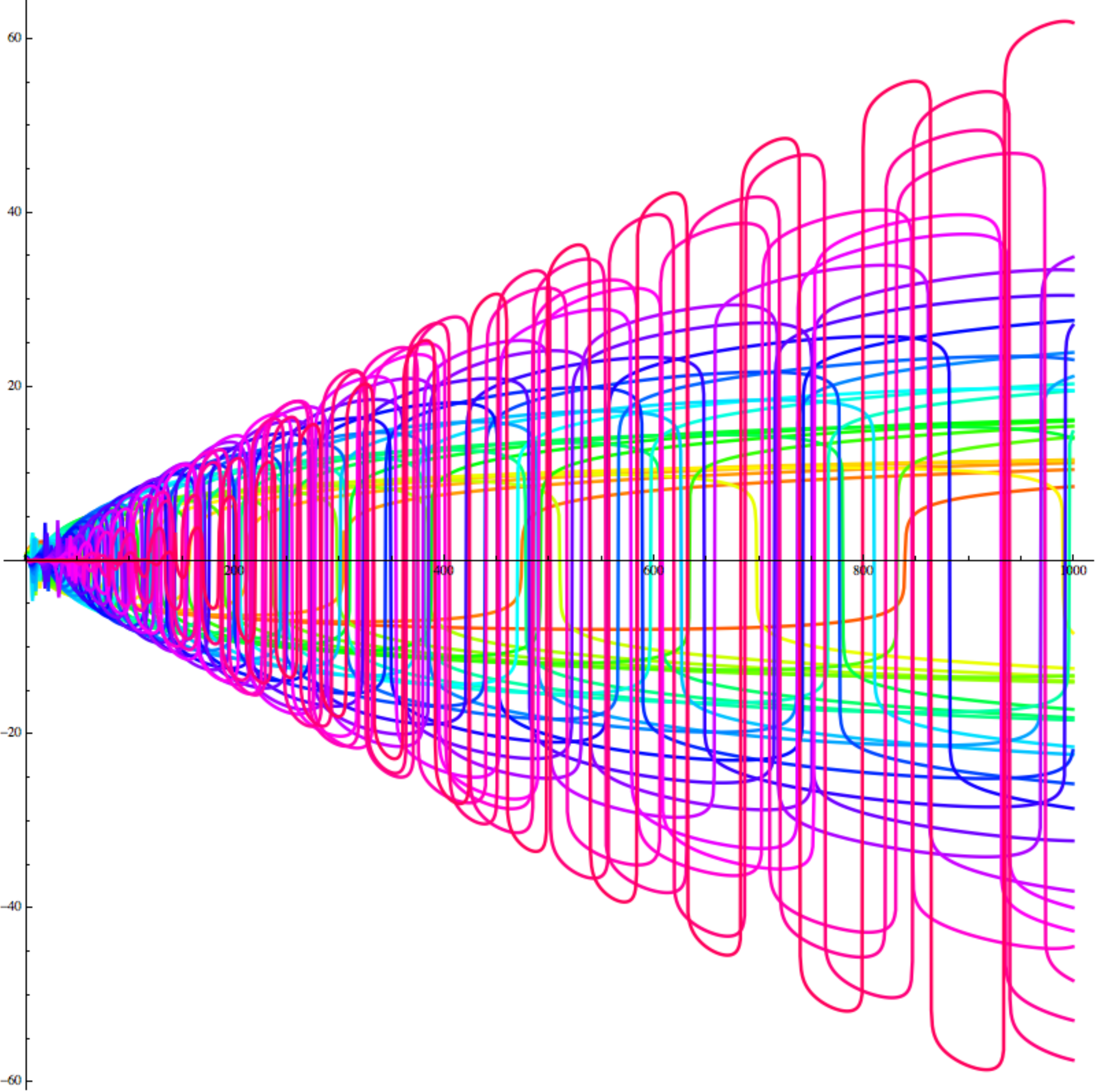}}}
}
\caption{
\label{figure2}
The expectation of the Euler characteristic $\chi$ is seen when plotting various functions
$n \to f_p(n) = \log^{\pm}({\rm E}_{n,p}[\chi])$ for $n \leq 100$ and $n \leq 1000$.}
\end{figure}

The simplest notion of a Ricci curvature $R(e)$ of an edge $e$ is the average over the curvatures of all wheel graph centers 
which contain $e$ as an edge connected to the center. The scalar curvature $S(v)$ at a vertex $v$ is then the average over 
all wheel graph curvatures which contain a vertex. 
For a connected graph, the vacuum Einstein equations $R-S=0$ are satisfied if and only if the wheel curvatures
are constant. As in the continuum, graphs with constant sectional curvature solve also the Einstein equations. \\

Do such graph have extremal $\chi(G)$? Looking at small $n$, it appears at first as if complete bipartite graphs 
$K_{n,n}$ with $\chi(K_{n,n}) = 2n-n^2$ could lead to the minimum among all graphs of order $2n$.
But this is not the case, as we can compute the expectation of the Euler characteristic on 
Erdoes-R\'enyi graphs explicitly in \cite{randomgraph} as 
${\rm E}_{n,p}[\chi] = \sum_{k=1}^n (-1)^{k+1} \B{n}{k} p^{\B{k}{2}}$. 
Figure \ref{figure2} shows a collection of functions
$n \to f_p(n) = \log^{\pm}({\rm E}_{n,p}[\chi])$ for $n \leq 100$ and $n \leq 1000$,
where $\log^{\pm}(x) = \sign(x) \log|x|$ and ${\rm E}_{n,p}$ is the expectation in the
Erdoes-Renyi probability space of graphs of order $n$ in which edges are turned on with
probability $p$.  In each case, we have plotted the function $f_p$ for fifty $p$ values
between $0$ and $1$ so that the hulls produce bounds for the extremal Euler characteristic.
The actual maxima or minima are outside the enclosed cone. \\

For $p=0.5$ and $n=300$ already, the Euler characteristic average has become
much smaller than for the bipartite graph $K_{n,n}$ and for $p=0.9$ and $n=400$ larger than 
$\chi(P_n)=n$ with no edges. The probabilistic argument is not constructive. But it shows that
the maximum Euler characteristic is larger than $e^{c n}$ and the minimal smaller than $-e^{c n}$ for some $c>0$.
For example, we are not able to give concrete graphs with $n=1000$ for which $\chi(G) > e^{50}$ 
even so we know that they exists as the expected value is larger. 
Computing the Euler characteristic of a large graph is a formidable task because if an edge is turned 
on with probability $p$ we already expect $p n(n-1)/2$ edges and even if we compute the 
Euler characteristic using the Poincar\'e-Hopf method \cite{poincarehopf}, a computer can not get an 
answer if $n=1000$ and say $p=1/2$. 

\begin{figure}[H]
\parbox{15.4cm}{
\parbox{7cm}{\scalebox{0.22}{\includegraphics{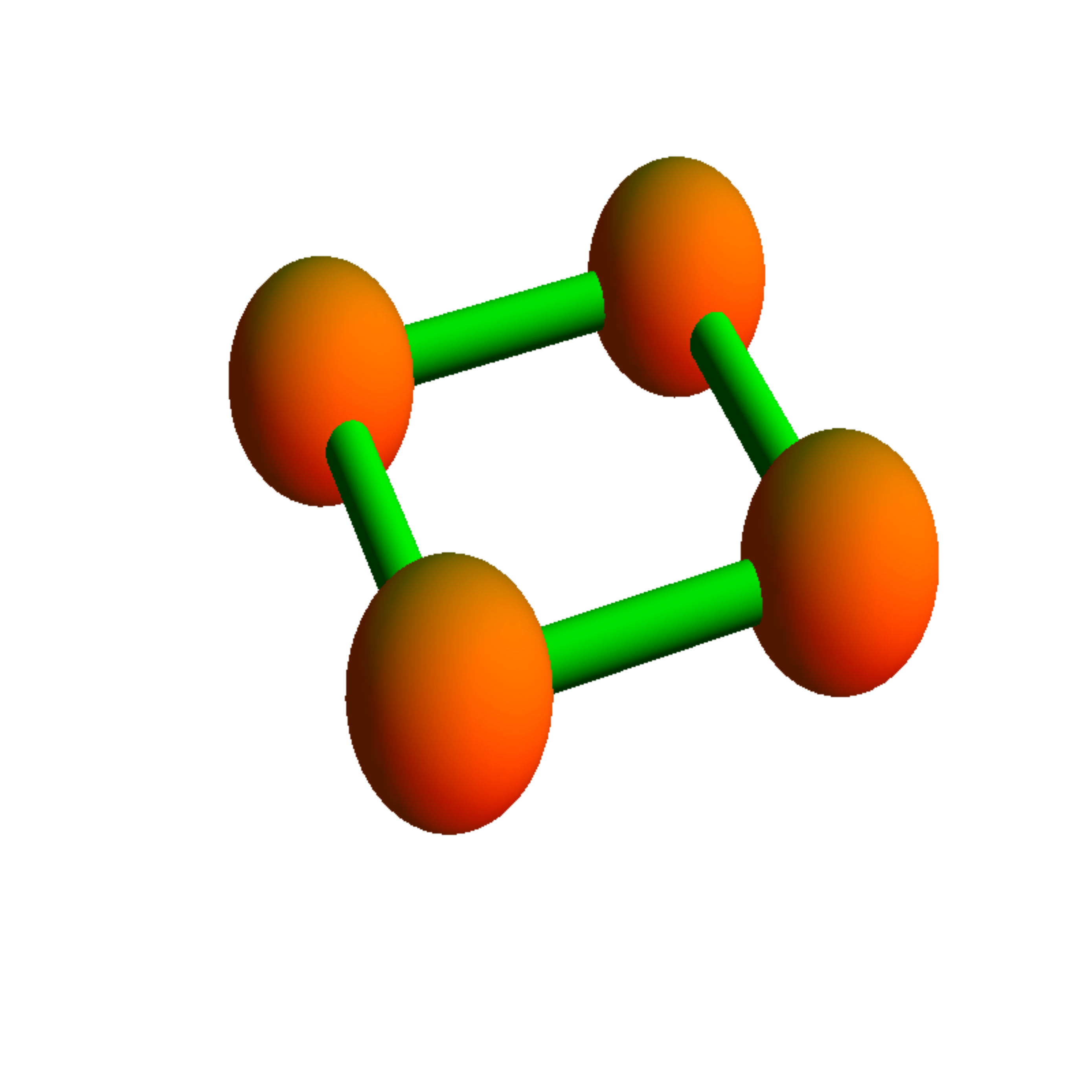}}}
\parbox{7cm}{\scalebox{0.22}{\includegraphics{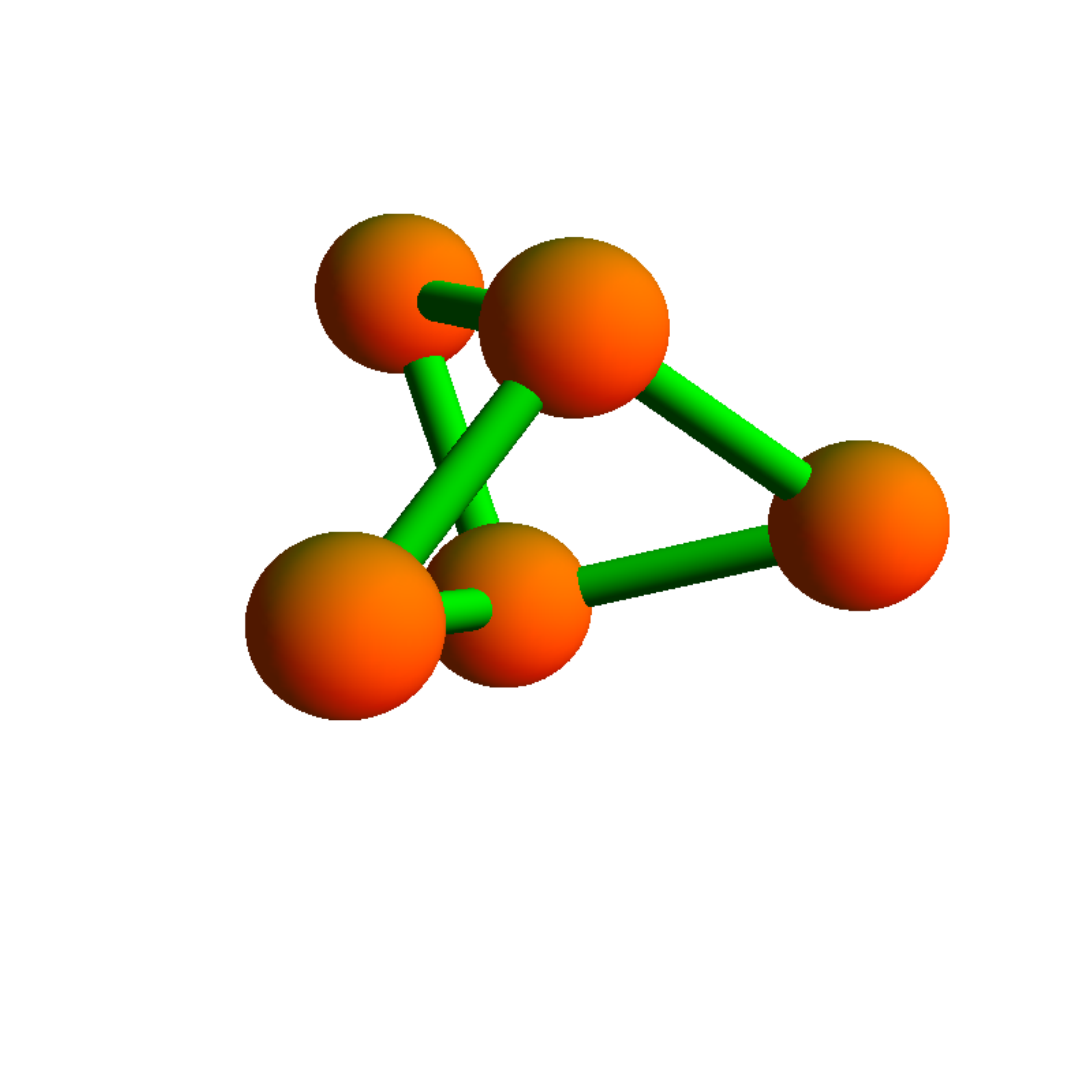}}}
}
\caption{
Two connected graphs with minimal Euler characteristic in the class of graphs of the same order.
$C_4$ is equal to $K_{2,2}$ and the two star graph $T_2$ with $\chi(T_2)=-1$ 
obtained by gluing two star graphs along the rays does not have constant curvature. The two
star-centers have curvature $=-1/2$. \label{figure1}
}
\end{figure}

Figures~(\ref{figure1}) and (\ref{figure2}) show examples of Sarumpaet graphs 
for small $n$. The hyperbolic graph with constant negative 
curvature $-1/2$ has minimal Euler characteristic $-3$ among all connected graphs of order $6$. 
More generally, the complete bipartite graph $K_{n,n}$ has no triangles and so $\chi(K_{n,n}) = 2n-n^2$.
The maximum for $n=6$ is the octahedron, which is a constant positive curvature $1/3$ graph.
The number $n=6$ is stable in the sense that $n=5$ both the minimum and maximum changes and for 
$n=7$, the minimum and maximum does not increase.
The minimum and maximum is monotone in $n$ because we can make homotopy deformations:
growing a single hair does not change the Euler characteristic.
Do geometric graphs $G$ with constant sectional curvature have extremal Euler characteristic? 
Torus graphs with zero curvature show that such graphs do not have to have maximal or minimal
Euler characteristic but extremal could mean ``stationary" in a more general sense, tori being
saddle type extrema. \\

\begin{figure}[H]
\parbox{15.4cm}{
\parbox{7cm}{\scalebox{0.22}{\includegraphics{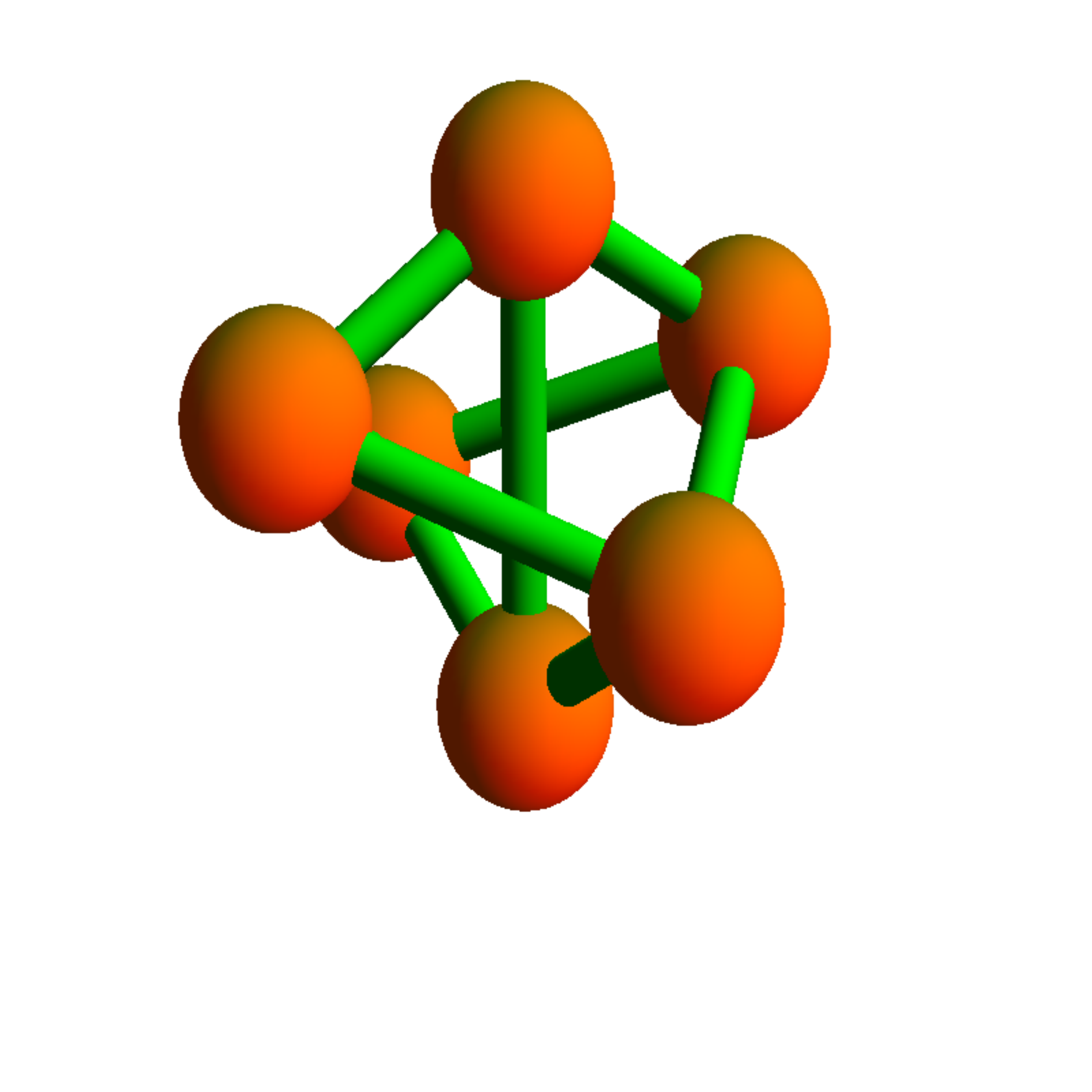}}}
\parbox{7cm}{\scalebox{0.22}{\includegraphics{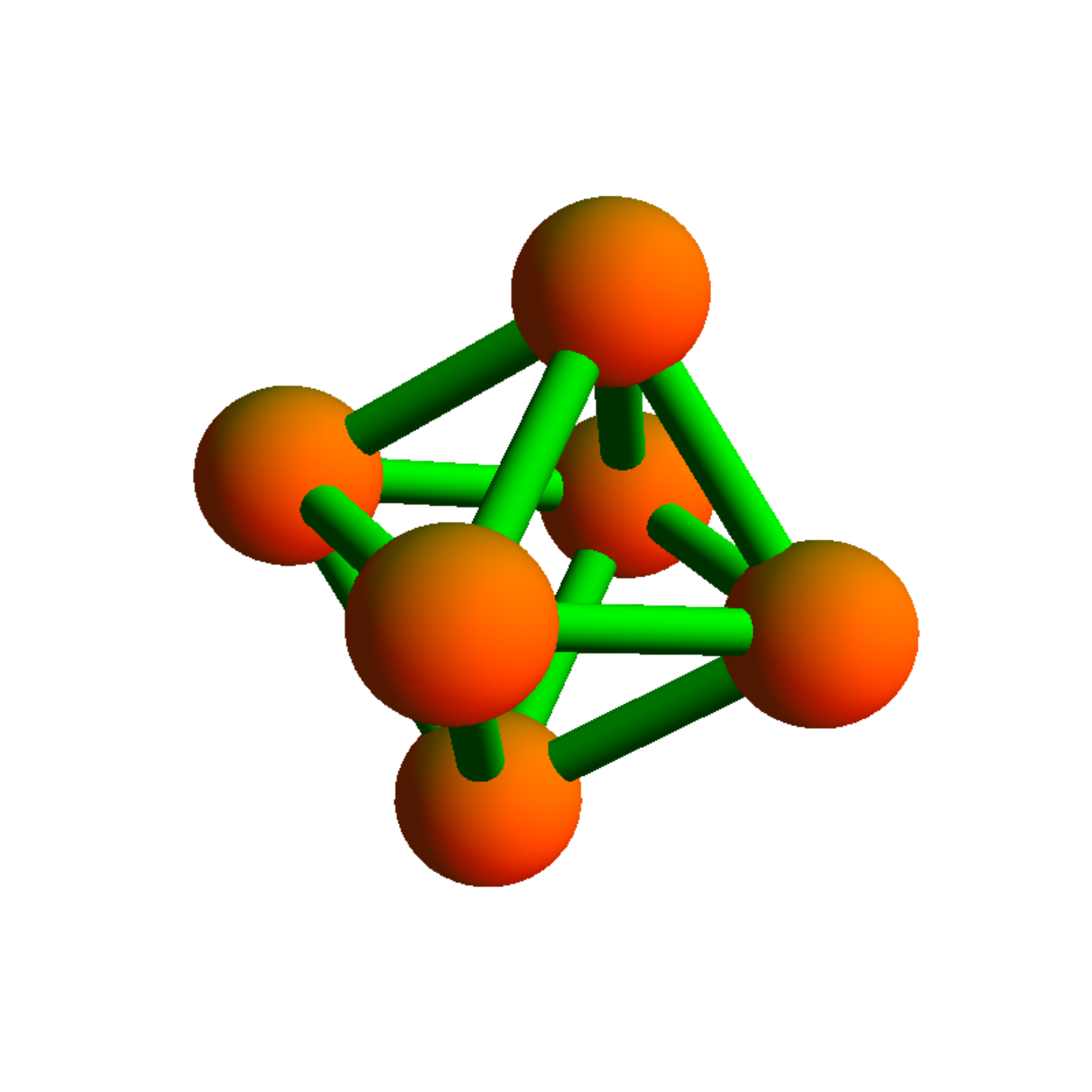}}}
}
\caption{
The hyperbolic utility graph $K_{3,3}$ is a minimum for $n=6$, the octahedron is a positive
curvature graph and is the maximum among all graphs of order $6$.
\label{figure2}
}
\end{figure}

More daring is to ask whether graphs with extremal Euler characteristic are geometric graphs (possibly 
up to trivial homotopies like diagonal flips) or whether certain dimensions are selected out. 
For any finite simple graph $G$, we have defined the Ricci curvature at an edge $e$ as the average
of curvatures of all wheel graphs containing $e$. The scalar curvature is a function on vertices and
averages all the Ricci curvatures of adjacent edges. The Einstein tensor $T_v(e)=R(e) - R(v)$ {\bf defines} then 
the mass tensor, a function on the edges attached to $x$. This is defined for any finite simple graph. 
By definition, the mass tensor satisfies the conservation law $\sum_{e \in B(v)} T_v(e) = 0$.
The metric tensor has not entered the above equations because the metric is still trivial. But letting the 
Noether symmetry group  \cite{IsospectralDirac2} act on the geometry, changes this. 
Since we know now to replace the scalar curvature by an average over curvatures of 
surfaces defined by a function $f$, we can redefine Ricci, scalar and mass tensor to have modified 
Einstein equations. The new ones are deformable under quantum deformations and still work for general 
finite simple graphs even so if the graphs are not symmetric, we also have to deal with the expectation
of the Euler characteristic of spheres. The new setup is more sensible to quantum mechanics or symmetries
which deform the metric: if the geometry changes through a unitary deformation of the Dirac operator, 
the Ricci curvature and mass move along nicely. It is a consequence of Gauss-Bonnet that the sum of
Ricci curvatures over all edges is related to $\chi(G)$ and that therefore the total mass 
satisfies a conservation law. \\

We started to investigate the Einstein equations for general finite simple graphs and look with 
the computer for Einstein graphs, graphs for which the trace-less Ricci curvature = Einstein tensor is zero.
Symmetric graphs like star graphs and circular graphs, regular polyhedra or complete graphs 
are Einstein. All the extrema mentioned here are Einstein but we see that also many Einstein graphs are 
not global maxima or minima. In which sense they can be seen as critical points still remains to be seen.


\end{document}